\documentclass[11pt,leqno]{amsart}
\usepackage{amsmath,epsfig}
\usepackage{amssymb}
\numberwithin{equation}{section}
\def\1{\ensuremath{\mathrm{1}\hspace{-.35em} \mathrm{1}}} 

\def\E{\mathbb{E}}

\def\N{\mathbb{N}}
\def\P{\mathbb{P}}
\def\R{\mathbb{R}}
\def\Z{\mathbb{Z}}
\def\vv{\mathop{\rm Var}\nolimits}
\def\cov{\mathop{\rm Cov}\nolimits}

\usepackage[T1]{fontenc}
\usepackage{amssymb}
\usepackage{dsfont}
\numberwithin{equation}{section}
\newtheorem{thm}{Theorem}[section]
\newtheorem{prop}[thm]{Proposition}
\newtheorem{lemma}[thm]{Lemma}

\theoremstyle{definition}
\newtheorem{dfn}{Definition}[section]

\theoremstyle{remark}

\allowdisplaybreaks

\begin{document}

\title[Deviation inequalities for  sums of weakly dependent time series]
{Deviation  inequalities for sums of weakly dependent time series}

\author[O. Wintenberger]{Olivier Wintenberger}

\address{CEREMADE, UMR CNRS 7534\\
Universit\'e de PARIS - DAUPHINE\\
Place du Mar\'echal De Lattre De Tassigny\\
75775 PARIS CEDEX 16 - FRANCE\\
{wintenberger@ceremade.dauphine.fr}}

\maketitle

\begin{abstract}
In this paper we give new deviation inequalities of Bernstein's type for the partial sums of weakly dependent time series.  The loss from the independent case is studied carefully. We give non mixing examples such that dynamical systems and Bernoulli shifts for whom our deviation inequalities hold. The proofs are based on the blocks technique and different coupling arguments. 
\end{abstract}


\section{Introduction}
The aim of this paper is to extend the deviation inequality of Bernstein's type from the independent case to some weakly dependent ones. We consider a sample $(X_1,\ldots,X_n)$ of a stationary process $(X_t)$ in a metric space $(\mathcal X,d)$. Considering the set $\mathcal F$ of  $1$-Lipschitz functions from $\mathcal X$ to $[-1/2,1/2]$, we are interested by the deviation of the partial sum $S(f)= \sum_{i=1}^n f(X_i)$ for any $f\in\mathcal F$ assuming that $\E(f(X_i))=0$. If the $X_i$ are independent and if $\sigma_k^2(f)=k^{-1}\vv(\sum_{i=1}^kf(X_i))$, the classical Bernstein inequality gives the deviation estimate, see Bennett \cite{Bennett1962}: 
\begin{equation}\label{deviid}
\P\Big(S(f)\ge \sqrt{2 n\sigma_1^2(f) x }+x/6\Big)\le e^{-x}\quad\mbox{for all}\quad x\ge0.
\end{equation}
This inequality reflects the gaussian approximation of the tail of $S(f)$ for small values $x$. And for large values of $x$, it reflects  the exponential  approximation of the tail of $S(f)$. This deviation inequality is very useful in statistics, see for example the monographs of Catoni \cite{Catoni2004} and of Massart \cite{Massart2006}.\\

To extend such deviation inequality to the dependent cases, a tradeoff between the sharpness of the estimates and the generality of the context has to be done. Estimates as sharp as in the independent cases (up to constants) are obtained for Markov chains in Lezaud \cite{Lezaud1998}, Joulin and Ollivier \cite{Joulin} under granularity. Bertail and Clemençon \cite{Bertail} obtain a deviation inequality for recurrent Markov chains. There exists $C>0$ such that  for all $M>0$ and all $x>0$:
$$
\P(S(f)\ge C(\sqrt{n\sigma^2(f)_Tx}+Mx))\le e^{-x}+n\P(T_1\ge M),
$$
where the $T_i$ are the iid regeneration times and $\sigma_T^2(f)=\E(T)^{-1}\vv(\sum_{i=1}^{T_1}f(X_i))$. Up to a constant, it is the limit variance in the CLT of $S(f)$, more natural than $\sigma^2_1(f)$ in \eqref{deviid}. This primitive estimate of the tail is natural as, through the splitting technique of Nummelin \cite{Nummelin1978}, the partial sums $S(f)$ are sums of iid sequences of blocks of size $T_i$. If the regeneration times are bounded, then up to different constants the same estimate than in the iid case is obtained. If the regeneration times  admit finite exponential moments, fixing $M\approx \ln n$ Adamczak \cite{Adamczak2008} obtains estimates of the deviations with a constant $C>0$:
$$
\P\Big(S(f)\ge C\Big(\sqrt{n \sigma^2_T(f) x }+\ln nx\Big)\Big)\le e^{-x}\quad\mbox{for all}\quad x\ge0.
$$
A loss of rate $\ln n$, that cannot be reduced, appears in the exponential approximation compared with the iid case, see Adamczak  \cite{Adamczak2008} for more details.\\

In all these works, the strong Markov property is crucial. To bypass the Markov assumption, one way is to use dependent coefficients. Ibragimov \cite{Ibragimov1962} introduced the uniformly $\phi$-mixing coefficients. In this settings, Samson \cite{Samson2000} achieves the deviation inequality \eqref{deviid} with different constants as soon as $\sum \sqrt\phi_r<\infty$.  Less accurate results have been obtained for more general mixing coefficients than the $\phi$-mixing ones: Viennet \cite{Viennet1997} for absolutely regular mixing and Merlevede {\it et al.} \cite{Merlevede} for geometrically strongly mixing. Recently, mixing coefficients have been extended to weakly dependent ones, see Doukhan and Louhichi \cite{Doukhan1999} and  Dedecker and Prieur \cite{Dedecker2005}. Under the exponential decrease of these coefficients, deviation inequalities for $S(f)$ with a loss in the exponential approximation are given in Doukhan and Neumann \cite{Doukhan2007}. Merlevede {\it et al.} \cite{Merlevedea} extends these results for the partial sum $S(f)$ for unbounded functions $f$.\\

The dependence context of this paper is the one of the so-called $\varphi$-weakly dependent coefficients introduced by Rio in \cite{Rio1996} to extend the uniformly $\phi$-mixing coefficients. We provide new deviation inequalities for non mixing processes, such that dynamical systems called expanding maps, see Collet et al. \cite{Collet2002} and continuous functions of Bernoulli shifts. The Bernstein's deviation inequality in these non mixing contexts sharpens the existing ones. The deviation inequality is obtained by dividing the sample $(X_1,\ldots,X_n)$ in different blocks $(X_i,\ldots,X_{i+k^\ast})$, where the length $k^\ast$ must be carefully chosen and then by approximating non consecutive blocks by independent blocks using a coupling scheme.\\

The coupling scheme follows from a conditional Kantorovitch-Rubinstein duality due to Dedecker {\it et al.} \cite{Dedecker2006} and detailed in Section \ref{prel}. Using this coupling argument,  a new deviation inequality can be stated in Section \ref{devin1}:
\begin{equation}\label{devinwd}
\P\Big(S(f)\ge 5.8\sqrt{n \overline\sigma_{k^\ast}^2(f) x }+1.5~k^\ast x)\Big)\Big)\le e^{-x}\quad\mbox{for all}\quad x\ge0,
\end{equation}
with $\overline \sigma_j^2(f) =\sup_{j\le k\le n}\sigma_k^2(f)$ for all $1\le j\le n$ and $k^\ast=\min\{ k\ge 1;~ k\delta_k\le \overline\sigma_k^2(f)\} $, where $(\delta_k)$ only depends of the $\varphi$-coefficients, see condition \eqref{MC} for more details. Unlike $ \sigma_1(f)$ in \eqref{deviid}, the variance term $\overline \sigma_{k^\ast}^2(f)$ is natural as it tends to the limit variance in the CLT with ${k^\ast}$. When the TLC holds, i.e. $ \overline\sigma_k^2(f)\to \sigma^2(f)>0$, then the classical Bernstein's inequality \eqref{deviid} holds up to constants with $\sigma_1^2(f)$ replaced by $\sigma^2(f)$, see Subsection \ref{remth} for more details. On the opposite, if the functionals $f_n$ are such that $\sigma_1^2(f_n)\to0$, then for exponentially decreasing rate of the $\varphi$-coefficients, $k^\ast\approx -\ln(\overline\sigma_1^2(f_n))$ and a logarithmic loss appears.  As in the recursive Markov chains case, the loss in the exponential approximation depends on the size of blocks $k^\ast$. We do not know if this loss in the exponential domain may be reduced for such non uniformly $\phi$-mixing sequences.\\

In many practical examples such that chains with infinite memory introduced by Doukhan and Wintenberger \cite{Doukhan2008}, Bernoulli shifts and Markov kernels, an $L^\infty$ coupling scheme is tractable, see Section \ref{lc} for a detailed definition. In these specific cases of $\varphi$-weakly dependent sequences, an improved version of the deviation inequality \eqref{devinwd} is given in Theorem \ref{thcl}:
$$
\P\left(S(f)\ge 2\sqrt{n\sigma^2_{{k^\ast}'}(f)x}+1.34~{k^\ast}' x)\right)\le e^{-x}\quad\mbox{for all}\quad x\ge0,
$$
with ${k^\ast}'=\min\{1\le k\le n~/~n\delta'_k\le kx\}$, where $(\delta_k')$ only depends of the $L^\infty$ coupling scheme, see condition \eqref{cl} for more details. The paper finishes with the proofs collected in Section \ref{pf}.

\section{Preliminaries: coupling and weak dependence coefficients}\label{prel}
Let $(X_1,\ldots,X_n)$ with $n\ge 1$  be  a sample of random variables on  some probability space $(\Omega,\mathcal A,\P)$ with value in a metric space $(\mathcal X,d)$. We assume in all the sequel that for any $n\ge 1$ there exists a strictly stationary process $(X_t^{(n)})$ such that $(X_1,\ldots,X_n)=(X_1^{(n)},\ldots,X_n^{(n)})$. Let us consider $\mathcal F$ the set of measurable functions $f:$ $\mathcal X\mapsto \R$ satisfying:
\begin{equation}\label{Lipb}
|f(x)-f(y)|\le d(x,y),\quad \forall (x,y)\in\mathcal X\times \mathcal X\quad\mbox{ and }\quad\sup_{x\in\mathcal X}|f(x)|\le 1/2.
\end{equation}
We denote the partial sum $S(f)=\sum_{i=1}^nf(X_i)$ and $\mathcal M_j=\sigma(X_t;1\le t\le j)$ for all $1\le j\le n$.

\subsection{Kantorovitch-Rubinstein duality}
The technique of coupling is related with the Kantorovitch-Rubinstein duality. The duality states that  given two distribution $P$ and $Q$ on $\mathcal X$ there exists a random couple $Y=(Y_1,Y _2)$ with $Y_1\sim P$ and $Y_2\sim Q$ satisfying
$$
\E (d(Y_1,Y_2))=\sup_{f\in\Lambda_1}\E|f(dP-dQ)|=\inf_{Y'}\E (d(Y'_1,Y'_2)),
$$
where $Y'$ have the same margins than $Y$ and $\Lambda_1$ denotes the set of $1$-Lipschitz functions such that $|f(x)-f(y)|\le d(x,y)$.\\

Dedecker, Prieur and Raynaud de Fitte \cite{Dedecker2006} extend the classicalKantorovitch-Rubinstein duality in the time series framework by considering it conditionally on some event $\mathcal M\in\mathcal A$. Assuming that the original space $\Omega$ is  rich enough, i.e. it exists a random variable $U$ uniformly distributed over $[0,1]$ and independent of $\mathcal M$, for any $Y_1\sim P$ with values in a Polish space it exists a random variable $Y_2\sim P$ independent of $\mathcal M$ satisfying
\begin{equation}\label{krcond}
\E(d(Y_1,Y_2)~|~\mathcal M)=\sup\{|\E(f(Y_1)|\mathcal{M})-\E(f(Y_1))|,
f\in\Lambda_1\}\qquad a.s.
\end{equation}

\subsection{$\varphi$-weak dependence coefficients and coupling schemes}
Let us recall the weak dependence coefficient $\varphi$ introduced in Rio \cite{Rio1996}
\begin{dfn}\label{varphi}
For any $X\in\mathcal{X}$, for any $\sigma$-algebra $\mathcal{M}$ of $\mathcal A$ then
$$
\varphi(\mathcal{M},X)=\sup\{\|\E(f(X)|\mathcal{M})-\E(f(X))\|_\infty,
f\in\mathcal F\}.
$$
\end{dfn}
Another equivalent definition is 
\begin{equation}\label{theta}
\varphi(\mathcal M,X_r))=
\sup \{|\cov(Y,f(X_r))|,f\in\mathcal F\mbox{ and }Y\mbox{ is }\mathcal M\mbox{-measurable and }\E |Y|=1\},
\end{equation}
see \cite{Dedecker2005}.\\

We will denote by {\bf (A)} the specific case where $\mathcal X$ is a Polish space with $\sup_{(x,y)\in\mathcal X^2}d(x,y)\le 1$. In the case {\bf (A)} we have  $\varphi(\mathcal M,X)=\tau_\infty(\mathcal{M},X)$ where $\tau_\infty$ is the coupling coefficient defined in \cite{Dedecker2007a} by the relation
$$
\tau_\infty(\mathcal M,X)=\sup\{\|\E(f(X)|\mathcal{M})-\E(f(X))\|_\infty,
f\in\Lambda_1\}.
$$
This last coefficients is the essential supremum of the right hand side term of the conditional Kantorovicth-Rubinstein duality \eqref{krcond}. Thus in the case  {\bf (A)} we get a coupling scheme directly on the variable $X$ via the Kantorovitch-Rubinstein duality \eqref{krcond}: it exists a version $X^\ast\sim X$ independent of $\mathcal M$ such that 
$$
\|\E(d(X,X^\ast)~|~\mathcal M)\|_\infty=\tau_\infty(\mathcal M,X)=\varphi(\mathcal M,X).
$$
When $d$ is the Hamming distance $d(x,y)=\1_{x\neq y}$  the coefficient $\varphi(\mathcal{M},X)$ coincides with the uniform mixing coefficient $\phi(\mathcal M,\sigma(X))$ of Ibragimov defined for 2 $\sigma$-algebra $\mathcal M$ and $\mathcal M'$ as:
$$
\phi(\mathcal M,\mathcal M')=\sup_{M\in\mathcal M,M'\in\mathcal M'}|\P(M'~|~M)-\P(M')|.
$$

In more general context than {\bf (A)}, we have $\varphi(\mathcal M,X)\le \tau_\infty(\mathcal{M},X)$ and coupling scheme directly on $X$ is not tractable. Then we do coupling scheme on the variables $f(X_i)$ for some function $f$: If the sample $(X_1,\ldots X_n)$ is such that the coefficients $\varphi(\mathcal M_j,X_i)$ are finite for $1\le j<i\le n$ and if $f\in\mathcal F$ then the coupling scheme for $f(X_i)$  follows from the conditional Kantorovitch-Rubinstein duality \eqref{krcond} and the relation
$$
\tau_\infty(\mathcal M_j,f(X_i))\le \varphi(\mathcal M_j,X_i):
$$
There exists $f(X_i)^\ast$ such that $f(X_i)^\ast\sim f(X_i)$ is independent of $\mathcal M_j$ and
$$
\|\E(|f(X_i)^\ast-f(X_i)|~|~\mathcal M_j)\|_\infty=\tau_\infty(\mathcal M_j,f(X_i)) \le \varphi(\mathcal M_j,X_i).
$$
In the case {\bf (A)} we also have another possible coupling scheme for $f(X_i)$, see Section \ref{lc} for practical examples: $f(X_i^\ast)\sim f(X_i)$ is independent of $\mathcal M_j$ and
$$
\|\E(|f(X_i^\ast)-f(X_i)|~|~\mathcal M_j)\|_\infty\le \|\E(d(X_i^\ast,X_i)|~|~\mathcal M_j)\|_\infty=\tau_\infty(\mathcal M_j,f(X_i)) = \varphi(\mathcal M_j,X_i).
$$

\subsection{Extensions on the product space $\mathcal X^q$, $q>1$.}

To consider conditional coupling schemes of length $q>1$ we need to extend the notions  of weak dependence coefficients on $X=(X_t)_{r\le t<r+q}\in\mathcal X^q$. It depends on the metric $d_q$ chosen for $\mathcal X^q$:
\begin{dfn}
For any $q\ge 1$, any $X\in \mathcal X^q$ and any $\sigma$-algebra $\mathcal M$ of $\mathcal A$ let us define the coefficients
$$
\varphi(\mathcal M,X)=\sup\{\|\E(f(X)|\mathcal{M})-\E(f(X))\|_\infty,
f\in\mathcal F_q\},
$$
where $\mathcal F_q$ is the set of $1$-Lipschitz functions with values in $[-1/2,1/2]$ of $\mathcal X^d$ equipped with the metric $d_q(x,y)=q^{-1}\sum_{ i=1}^qd(x_i,y_i)$. 
\end{dfn}
Let us discuss the consequences of the choice of the metric $d_q$:
\begin{itemize}
\item The $\tau_\infty$ coupling coefficients on $\mathcal X^q$ are defined for the metric $d_q$, see \cite{Dedecker2007a}, and for all $f\in\mathcal F$:
$$
\tau_\infty(\mathcal M,(f(X_1),\ldots,f(X_q)))\le  \varphi(\mathcal M,(X_1,\ldots,X_q)).
$$
Moreover, in the case {\bf (A)} it holds $\tau_\infty(\mathcal M,X)= \varphi(\mathcal M,X)$.
\item If $d$ is the Hamming metric, as $d_q(x,y)\le \1_{x\neq y}$ then $\varphi(\mathcal M,X)\le\phi(\mathcal M,\sigma(X))$. Thus the definition of the weakly dependent coefficients $\varphi$ differs here from the one of Rio in \cite{Rio1996} where $\mathcal X^q$ is equipped with $d_\infty(x,y)=\max_{1\le i\le q}d(x_i,y_i)$.
\end{itemize}

\subsection{First application: deviation inequality of Hoeffding type}

This application is due to Dedecker and Prieur  \cite{Dedecker2005}. Assume that $\mathcal X$ is a Polish space such that $\sup_{x,y}d(x,y)\le 1$, i.e. we are  in the case {\bf (A)}.
Assume that the coefficients $\varphi(\mathcal M_j,(X_{j+1},\ldots,X_n))$ are finite for all $1\le  j\le n-1$ and that $g:\mathcal X^n\to \R$ satisfies
$$
|g(x_1,\ldots,x_n)-g(y_1,\ldots,y_n)|\le \sum_{i=1}^nd(x_i,y_i).
$$
As we are  in the case {\bf (A)} it exists a coupling scheme $(X_{j+1}^\ast,\ldots,X_n^\ast)$ of $(X_{j+1},\ldots,X_n)$ for any $1\le j\le n-1$ such that, keeping the same notation than in \cite{Rio2000a}:
\begin{multline*}
\Gamma(g)=\|\E(g(X_{j+1},\ldots,X_n)~|~\mathcal M_j)-\E(g(X_{j+1},\ldots,X_n)\|_\infty\\
=\|\E(g(X_{j+1},\ldots,X_n)-g(X_{j+1}^\ast,\ldots,X_n^\ast)~|~\mathcal M_j)\|_\infty
\le (n-j)\varphi(\mathcal M_j,(X_{j+1},\ldots,X_n)).
\end{multline*}
Applying Theorem 1 of \cite{Rio2000a}, if $\E(g(X_1,\ldots,X_n))=0$ then for all $x\ge 0$ it holds:
$$
\P\left(f(X_1,\ldots,X_n)\ge \sqrt{2^{-1}\sum_{j=1}^n\left(1+2(n-j)\varphi(\mathcal M_{j},(X_{j+1},\ldots,X_n))\right)^2x}\right)\le e^{-x}.
$$
This deviation inequality of Hoeffding type only differs from the one for independence  by a constant. However, such inequalities are not as satisfactory as Bernstein ones for statistical applications.

\section{Deviation inequality around the mean inequality}\label{devin1}

Let us give an inequality for the deviation around the mean of $S(f)=\sum_{i=1}^nf(X_i)$ for $f\in\mathcal F$, with $(X_1,\ldots,X_n)$ on the metric space $(\mathcal X,d)$ and such that there exists a non increasing sequence $(\delta_r)$ that satisfies 
\begin{equation}
\label{MC}\sup_{1\le j\le n-2r+1}\varphi(\mathcal M_j,,(X_{r+j},\ldots,X_{2r+j-1}))\le \delta_r \mbox{ for all }r\ge 1.
\end{equation}

\subsection{A deviation inequality of Bernstein type}

Assume with no loss of generality that $\E(f(X_1))=0$.
\begin{thm}\label{mainsum}
For any  integer $n$, if there exists $(\delta_r)$ as in \eqref{MC}  then 
$$
\P\left(S(f)\ge  5.8\sqrt{n\overline \sigma_{k^\ast}^2(f) x}+1.5~k^\ast x\right)\le e^{-x},
$$
where
$
k^\ast=\min\{1\le k\le n~/~k\delta_k\le\overline\sigma^2_k(f)\}$ and  $\overline \sigma_{k^\ast}^2(f)=\max\{\sigma_{k}^2(f)~/~k^\ast\le k\le n \}$.
\end{thm}
The proof of this Theorem is given in Subsection \ref{pfmainsum}. We adopt the convention $\min \emptyset =+\infty$ and the estimate is non trivial when $r\delta_r\to 0$ and $n\delta_n\ge \overline \sigma_n(f)$, i.e. for not too small values of $n$. \\
 
Remark that the variance term $\overline \sigma_{k^\ast}^2(f)$ is more natural than $\sigma_1^2(f)$ in \eqref{deviid} as in the central limit theorem $\overline \sigma_{k^\ast}^2(f)$ converges to the limit variance as $k^{\ast}$ goes to infinity. Before giving some remarks on this Theorem, the next proposition give estimates of the quantity $\sigma^2_k(f)=k^{-1}\vv\left(\sum_{i=1}^kf(X_i)\right)$. 

\subsection{The variance terms $\sigma^2_k(f)$}

Under suitable assumptions on $(\delta_r)$, it is always possible to obtain rough estimates of $ \sigma_k^2(f)$ in function of $\sigma_1^2(f)$ and $\E|f(X_1)|$:
\begin{prop}\label{varest}
If the condtion \eqref{MC} is satisfied then we have for all $k\le n$ the inequality:
$$
\sigma^2_k(f)\le \left(\sigma_1^2(f)+2\E|f(X_1)|\sum_{r=1}^{k-1}\delta_r\right).
$$
\end{prop}
See Subsection for a straightforward proof \ref{pfvarest} of this Proposition. The estimate given in Proposition \ref{varest} can be rough, for example in the degenerate cases when $\sigma^2_k(f)$ tends to $0$ with $k$.

\subsection{Remarks on Theorem \ref{mainsum}}\label{remth}

The gaussian behavior around the mean is, up to a universal constant, the same than in the iid case with the more natural variance $\overline \sigma^2_{k^\ast}(f)$ instead of $\sigma_1^2(f)$ in \eqref{deviid}. However in the exponential domain the estimates given in Theorem \ref{mainsum} is sometimes less sharp than the one obtained for $\phi$-mixing in Samson \cite{Samson2000}. \\

 In the non degenerate case $\sigma^2_k(f)\to \sigma^2(f)>0$ then $k^\ast$ is finite  as soon as $r\delta_r\downarrow 0$. The deviation inequality of Theorem \ref{mainsum} becomes similar than the one in the iid case \eqref{deviid} with the  variance term $ \sigma^2(f)$ instead of $ \sigma^2_1(f)$: there exists $C>0$ such that for $n$ sufficiently large we have
$$
\P(S(f)\ge C(\sqrt{n\sigma^2(f)x}+x))\le e^{-x}\mbox{for all }x\ge 0.
$$
However, the estimate of the exponential behavior in Theorem \ref{mainsum} may differ from the one of the iid case. For example, for statistical issues it is often assumed that $f$ is chosen depending on $n$ such that $\sigma_1^2(f_n)\to 0$. Assume that $r\delta_r$ is summable. Using Proposition \ref{varest} and Jensen's inequality we have the estimate $\overline\sigma_1^2(f_n)\lesssim\sigma_1^2(f_n)^{1/2}$. If $\sigma_1^2(f_n)^{-1/2}n\delta_n\downarrow 0$ then for $n$ sufficiently large such that $k^\ast_n=\min\{k\le n~/~k\delta_k\le \sigma_1^2(f_n)^{1/2}\}$ exists, it holds
$$
\P(S(f)\ge C(\sqrt{\overline\sigma_{k^\ast_n}^2(f_n)nx}+k^\ast_n x))\le e^{-x}\mbox{ for all }x\ge 0\mbox{, with }C>0.
$$
As $k_n^\ast\uparrow\infty$ there is a loss compare with the iid case \eqref{deviid}. We do not know if this loss may be reduced outside the cases of uniformly mixing processes where \eqref{deviid} holds, see Samson \cite{Samson2000}.\\

This loss may be reduced when the autocorrelations are controlled, choosing a smaller size of blocks $k^\ast_n$. Assume that $\overline\sigma_1^2(f_n)\lesssim\sigma_1^2(f_n)$  (such relation is satisfied in the uniformly $\phi$-mixing context). If  $\sigma_1^2(f_n)^{-1}n\delta_n\downarrow 0$ then for $n$ sufficiently large such that $k^\ast_n=\min\{k\le n~/~k\delta_k\le \sigma_1^2(f_n)\}$ exists, it holds
$$
\P\left(S(f_n)\ge C\left(\sqrt{\sigma_1^2(f_n)nx}+k^\ast_n x\right)\right)\le e^{-x}\mbox{ for all }x\ge 0\mbox{, with }C>0.
$$
The loss compare with the iid case is due to $k^\ast_n\uparrow \infty$. More precisely
\begin{itemize}
\item If $\delta_r= C\delta^r$ for $C>0$ and $0<\delta<1$ then $k^\ast_n\approx -\ln(\sigma_1^2(f_n))$,
\item If $\delta_r= Cr^\delta$ for $C>0$ and $\delta>1$ then $k^\ast_n\approx \sigma_1^2(f_n)^{1/(1-\delta)}$.
\end{itemize}

\section{Examples}\label{ex}

We focus on non $\phi$-mixing examples as for them the inequality \eqref{deviid} holds up to constants, see Samson \cite{Samson2000}. We present dynamical systems that are known to be non $\phi$-mixing processes but they satisfy \eqref{MC} in the case {\bf (A)}. Other examples in the case {\bf (A)} are presented in the Section \ref{lc} as a sharpened deviation inequality holds for them, see Theorem \ref{thcl}. We also present in this Section continuous functions of Bernoulli shifts that are examples not $\phi$-mixing and not in the case {\bf (A)} and thus cannot be treated by the approach of Section  \ref{lc} and of \cite{Samson2000}.

\subsection{Dynamical systems}

Here we are in the case $\mathcal X=[0;1]$ and $d(x,y)=|x-y|$, i.e. in the {\bf (A)} and then $\varphi=\tau_\infty$. Since Andrews \cite{Andrews1984}, dynamical systems, defined as stationary solutions of $X_t=T(X_{t+1})$ for all $t$ are classical examples of non-mixing processes. Let us consider $X_t$ the stationary solution of
$$
X_t=\frac{1}{2}(X_{t-1}+\xi_t)
$$
where $(\xi_t)$ is an iid sequence distributed as a Bernoulli($1/2$). Then $X_t=T(X_{t+1}$ where $T(x)=2x$ modulo $1$. Even if it is not mixing, easy computation shows that $(X_t)$ satisfies \eqref{MC} with $r\delta_r= (4/9) 2^{-r}$ (in fact this specific case satisfies also  $r\delta_r'= (4/9) 2^{-r}$, see Section \ref{lc} for more details).\\

More general examples of dynamical systems are studied in Collet {\it et al.} \cite{Collet2002}. They obtain estimates of covariances terms, multivariate versions of \eqref{theta}, for dynamical systems called expanding maps. Then it follows the existence of $C>0$ and $0<\rho<1$ such that \eqref{MC} is satisfied with $r\delta_r= C \rho^{r}$, see Dedecker and Prieur \cite{Dedecker2005} for more details.

\subsection{Continuous functions of Bernoulli shifts}

Let us consider a $\phi$-mixing stationary process $(\xi_t)$ in some measurable space $\mathcal Y$ and a sequence $(U_t)$ in the metric space $\mathcal X$ defined as
$$
U_t=F(\xi_{t-j}; ~j\in\N),
$$
where $F$ is a measurable function. Assume that the original state space is large enough such that it exists $(\xi_t')$ distributed as $(\xi_t)$ but independent of it. As in \cite{Rio1996}, assume that there exists a non increasing sequence $(v_k)$ satisfying almost surely
$$
d(F(\xi_j;j\in\N),F(\xi^k_j;j\in\N))\le v_k,
$$
with the sequence $(\xi^k_t)$ satisfying $\xi_t=\xi_t^k$ for $0\le t\le k$ and for $t>k$, $\xi_t^k=\xi_t'$. Finally set $X_t=H(U_t)$ for some measurable function $H:\mathcal X\to \mathcal X$ and $t=\{1,\ldots,n\}$ and denote 
$$
w_H(x,\eta)=\sup_{d(x,y)\le \eta}d(H(x),H(y)).
$$
\begin{prop}\label{cb}
The sample $(X_1,\ldots,X_n)$ satisfies \eqref{MC} with
$$
\delta_r=\inf_{1\le k\le r-1}\{2\phi_{r-k}+\E(3w_H(U_0,2v_k))\wedge1\}.
$$
\end{prop}
See the Subsection \ref{pfcb} for the proof of this Proposition. Remark that by construction the process $(X_t)$ is non necessarily in the case ${\bf (A)}$

\section{In the case {\bf (A)} with a coupling scheme in $L^\infty$.}\label{lc}
In all this section we place us in the case {\bf (A)} where $\mathcal X$ is a Polish metric space and $d(x,y)\le 1$ for all $x,y\in\mathcal X$. For all $r\ge 1$ a coupling scheme in $L^\infty$ for $(X_i)_{r+j\le i<2r+j-1}$, $j\ge 1$, exists when we can construct $(X_i^\ast)_{r+j\le i<2r+j-1}$ distributed as $(X_i)_{r+j\le i<2r+j-1}$ and independent of $\mathcal M_j$ such that
\begin{equation}\label{cl}
\sup_{1\le j\le n-2r+1}\sum_{i=r+j}^{2r+j-1}d(X_i,X_i^\ast)\le r\delta_r'\quad a.s.\quad\mbox{for all }r\ge 1.
\end{equation}

\subsection{A sharper deviation inequality of Bernstein's type}

Remark that condition \eqref{cl} with $(\delta_r')$ implies condition \eqref{MC} with $\delta_r=\delta_r'$.
Then we obtain a slightly sharper deviation inequality than in Theorem \ref{mainsum}:
\begin{thm}\label{thcl}
For $f\in\mathcal F$ such that $\E(f(X_1))=0$ then we have for all $x\ge n\delta_k'$ and all $1\le k\le n$:
$$
\P(S(f)\ge x)\le \exp\left(-\frac{2n\sigma_k^2(f)}{k^2}h\left(\frac{k(x-n\delta_k')}{2n\sigma^2_k(f)}\right)\right)
$$
where $h(x)=(1+x)\ln(1+x)-x$ for all $x\ge 0$. Then it holds for all $x\ge 0$:
$$
\P\left(S(f)\ge 2\sqrt{n\sigma^2_{{k^\ast}'}(f)x}+1.34~{k^\ast}' x)\right)\le \exp(-x)
$$
with ${k^\ast}'=\min\{1\le k\le n~/~n\delta'_k\le kx\}$.
\end{thm}
The proof of this Theorem is given in Subsection \ref{pfthcl}.\\

Let us compare this deviation inequality with the result of Theorem \ref{mainsum}. In Theorem \ref{thcl} the variance term $\sigma^2_{k}(f)$ sharpens $\overline \sigma^2_{k}(f)$ and the normal approximation is better here. For the exponential approximation, in both Theorems losses are due to the chosen blocks sizes. As ${k^\ast}'=\min\{1\le k\le n~/~k\delta'_k\le xk^2/n\}$, if $k\delta_k'$ is decreasing as $k\delta_k$ then ${k^\ast} '\le k^\ast$ as soon as $n\sigma_{k}^2(f)\le k^2x$ or equivalently  $\sqrt{n\sigma_{k}^2(f)x}\le kx$, i.e. as soon as $x$ is in the domain of the exponential approximation. Thus for the normal and the exponential approximations, he deviation inequality in Theorem \ref{thcl} improves the one of Theorem \ref{mainsum}.\\ 

A tradeoff between the generality of the context and the sharpness of the deviation inequalities is done. Even if  \eqref{cl} is less general than \eqref{MC}, it is satisfied for many examples, see below.

\subsection{Bounded Markov Chains}
Following Dedecker and Prieur \cite{Dedecker2005}, let us consider a stationary Markov chain $(X_t)$ with transition kernel $P$ satisfying, for all $f\in\Lambda_1$, that $P(f)=\int f(y)P(x,dy)$ is a $\kappa$-Lipschitz function with $\kappa<1$. Then 
$$
r\delta_r'= \kappa^r(1+\cdots+\kappa^r),
$$
see  \cite{Dedecker2005} for more details. 

\subsection{Bounded chains with infinite memory}\label{chaininf}
 Let the sequence of the innovations $(\xi_{t})_{t\in{Z}}$ be an iid process on a measurable space $\mathcal{Y}$. We define
$X=(X_{t})_{t\in{Z}}$ as the solution of the equation
\begin{equation}\label{eq::rec}
X_t=F(X_{t-1},X_{t-2},\ldots;\xi_t)\qquad a.s.,
\end{equation}
for some bounded function $F:\mathcal X^{({N}\setminus\{0\})}\times \mathcal Y\to \mathcal{X}$ satisfying the condition
\begin{equation}\label{condlip}
d(F((x_k)_{k\in{N}\setminus\{0\}};\xi_0),F((y_k)_{k\in{N}\setminus\{0\}};\xi_0))
\le\sum_{j=1}^\infty a_j(F) d(x_j,y_j),~a.s.
\end{equation}
 for all
$(x_k)_{k\in{N}\setminus\{0\}}$, $(y_k)_{k\in{N}\setminus\{0\}}\in \mathcal{X}^{{N}\setminus\{0\}}$ such that there exists $N>0$ as $x_k=y_k=0$ for all $k>N$ and with $a_j(F)\ge 0$ satisfying
\begin{equation}
\label{condcontract} \sum_{j=1}^{\infty}a_j(F)
:=a(F) <1.
\end{equation}
Let $(\xi^\ast_t)_{t\in{Z}}$ be a stationary sequence  distributed
as $(\xi_t)_{t \in {{Z}}}$, independent of $(\xi_t)_{t \leq 0}$ and such that $\xi_t=\xi^\ast_t$ for $t>0$.
Let $(X^\ast_t)_{t\in\Z}$ be the solution of the equation
$$
X^\ast_t=F(X^\ast_{t-1},X^\ast_{t-2},\ldots;\xi_t^\ast),\qquad a.s.
$$
Using similar arguments than in Doukhan and Wintenberger \cite{Doukhan2008} we have the following result,
\begin{lemma}
\label{lemma_ex_1}
Under condition \eqref{condcontract} there exists some
bounded (by 1/2) stationary process $X$ solution of the equation \eqref{eq::rec}. Moreover, this solution satisfies \eqref{cl} with
$$r\delta_r'= \sum_{j=r}^{2r-1}\inf_{0<p\le j}\left\{a(F)^{r/p}+
\sum_{j=p}^\infty a_j(F)\right\}. $$
\end{lemma}
As the proof of this Lemma is similar than the one in \cite{Doukhan2008}, it is omitted here.\\
 
Many solutions of econometrical models may be written as chains with infinite memory. However, the assumption of boundedness is very restrictive for practical models.

\subsection{Bernoulli shifts}
Solutions of the recurrence equation \eqref{eq::rec} may always be written as $X_t=H((\xi_j)_{j\le t})$ for some measurable function $H:\mathcal Y^{\N}\mapsto \mathcal X$ were $(\xi_t)$ is an iid process called the innovations. In this very general framework, a coupling version $X_t^\ast$ is given by $X_t^\ast=H((\xi^\ast_t))$ where $(\xi^\ast_t)$ is a stationary sequence  distributed
as $(\xi_t)$, independent of $(\xi_t)_{t \leq 0}$ and such that $\xi_t=\xi^\ast_t$ for $t>0$. If there exist $a_i\ge 0$ such that
$$
d(H(x),H(y))\le \sum_{i\ge 1}a_id(x_i,y_i)\quad \mbox{with}\quad \sum_{i\ge 1}a_i<\infty,
$$
and if $\mathcal Y$ is a metric space such that it exists $y\in\mathcal Y$ with $d(\xi_1,y)$ bounded a.s., then $(X_t)$ satisfies \eqref{cl} with
$$
r\delta_r'=C\sum_{i\ge r}a_r
$$ 
for some $C>0$.

\section{Proofs}\label{pf}
This Section contains the proofs. 

\subsection{Proofs of the Theorems \ref{mainsum}}\label{pfmainsum} 
This section contains the proofs of the Bernstein's type estimates on the partial sums $S(f)$ for $f\in\mathcal F$. As in the independent case, the proofs follow the Chernoff device. 
We will proceed using Bernstein's block technique as in \cite{Doukhan1994}. Let us denote by $I_j$ the $j$-th block of  indices of size $k$, i.e.  $\{(j-1)k+1,jk\}$ except the last blocks and let $p$ be an integer such that $2p-1\le k^{-1}n\le 2p$.\\

Let us denote by $S_1$ and $S_2$ the sums of even and odd blocks defined as
$$
S_1=\sum_{i\in I_{2j},\, 1\le j\le p}f(X_i)\quad\mbox{ and }\quad S_2=\sum_{i\in I_{2j-1},\, 1\le j\le p}f(X_i).
$$
From Cauchy-Schwartz inequality, it holds:
$$
\ln \E [\exp(tS(f))]\le \frac{1}{2}\left(\ln \E \exp\left(2tS_1\right)+\ln \E \exp\left(2tS_2\right)\right).
$$
Now let us treat in detail the term depending on $S_1$, the same argument applies identically to  $S_2$. We want to prove that for any $0\le t\le 1$, choosing $k=[1/t]\wedge n$ as in \cite{Doukhan1994} it holds:
\begin{equation}\label{ineg1}
\ln \E(\exp(tS(f)))\le 4nt^2(2(e-2)\sigma^2_k(f)+ek\delta_k).
\end{equation}
Denoting $L_m=\ln \E(\exp(2t\sum_{i\in I_{2j},\, 1\le j\le m}f(X_i)))$ for any $1\le m\le p$, we do a recurrence on $m$ remarking that $\ln \E(\exp(2tS_1))=L_p$. From Holder inequality, we have for any $2\le m\le p-1$ the inequalities:
\begin{eqnarray*}
\exp(L_{m+1})-\exp(L_{m})
\exp(L_{1})&&\\ &&\hspace{-4cm}\le \exp(L_{m})\left\|\E\Big(\exp\Big(2t\sum_{i\in I_{2(m+1)}}f(X_i)\Big)~|~\mathcal M_{2mk}\Big)-\E\Big(\exp\Big(2t\sum_{i\in I_{2(m+1)}}f(X_i)\Big)\Big)\right\|_\infty\\
&&\hspace{-4cm} \le \exp(L_{m})\left\|\E\Big(\exp\Big(2t\sum_{i\in I_{2(m+1)}}f(X_i)\Big)-\exp\Big(2t\sum_{i\in I_{2(m+1)}}f(X_i)\Big)^\ast~|~\mathcal M_{2mk}\Big)\right\|_\infty,
\end{eqnarray*}
where $\exp\Big(2t\sum_{i\in I_{2(m+1)}}f(X_i)\Big)^\ast$ is a coupling version of the variable $\exp\Big(2t\sum_{i\in I_{2(m+1)}}f(X_i)\Big)$, independent of $\mathcal M_{2mk}$. From the definition of the coupling coefficients $\tau_\infty$, we know that
\begin{multline*}
\left\|\E\Big(\exp\Big(2t\sum_{i\in I_{2(m+1)}}f(X_i)\Big)-\exp\Big(2t\sum_{i\in I_{2(m+1)}}f(X_i)\Big)^\ast~|~\mathcal M_{2mk}\Big)\right\|_\infty\\
\le \tau_\infty\Big(\mathcal M_{2mk},\exp\Big(2t\sum_{i\in I_{2(m+1)}}f(X_i)\Big)\Big).
\end{multline*}
As $\sum_{i\in I_{2(m+1)}}f(X_i)$ is bounded with $k/2$, then $u\to \exp(2tu)$ is a Lipschitz function with constant $2kt\exp(kt)$ with respect to $d_k$ and bounded with $\exp(kt)$ for all $t\ge 0$. We then deduce that for $n^{-1}< t\le 1$, choosing $k=[1/t]\wedge( n-1)$ and under condition \eqref{MC} we have
$$
\tau_\infty\Big(\mathcal M_{2mk},\exp\Big(2t\sum_{i\in I_{2(m+1)}}f(X_i)\Big)\Big)\le  2kte^{kt}\varphi(\mathcal M_{2mk},(X_i)_{i\in I_{2(m+1)}})
\le 2e\delta_k.
$$
Collecting this inequalities, we achieve that
$$
\exp(L_{m+1})\le\exp(L_m)(\exp(L_1)+2e\delta_k). 
$$
The classical Bennett's inequality on $\sum_{i\in I_2}f(X_i)$ gives the estimates $\exp(L_1)\le 1+ 4\sigma^2_k(f)/k(e^{kt}-kt-1)$ and as $kt\le 1$ we obtain
$$
L_{m+1}\le L_m+\ln\left(1+\frac{4(e-2)\sigma^2_k(f)+2ek\delta_k}{k}\right)
\le L_m+ \frac{4(e-2)\sigma^2_k(f)+2ek\delta_k}{k}.
$$
The $p$ steps of the recurrence leads to the desired inequality
$$
\ln \E(\exp(2tS_1))\le 2p\frac{2(e-2)\sigma^2_k(f)+ek\delta_k}{k}.
$$
As the same inequality holds for $S_2$ we obtain \eqref{ineg1} for $n^{-1}<t\le 1$ remarking that $2pk^{-1}\le 4nt^2$. For $t\le n^{-1}$, classical Bennett inequality on $S_1$ gives
$$
\ln \E(\exp(2tS_1))\le 4\sigma_n^2(f)/n(e^{nt}-nt-1).
$$
Remarking that $e^{nt}-nt-1\le (nt)^2\sum_{k\ge0}(nt)^k/(k+2)!$ and $(k+2)!\ge 23^{k}$ we derive that $e^{nt}-nt-1\le 2^{-1}(nt)^2\sum_{k\ge0}3^{-k}\le 3/4(nt)^2$ for $nt\le 1$. Then collecting thes bounds, for $t\le n^{-1}$ it holds
$$
\ln \E(\exp(2tS_1))\le 3n\overline \sigma_{n}^2(f)t^2\le 4nt^2(2(e-2)\sigma^2_n(f)+en\delta_n).
$$
The same holds for $S_2$ and then \eqref{ineg1} follows for $0\le t\le n^{-1}$ and then for all $0\le t\le 1$.
\\

Note that for $k\ge k^\ast$ we have $\sigma^2_k(f)\le \overline  \sigma^2_{k^\ast}(f)$ and $k\delta_k\le \overline\sigma^2_k(f)$ by definition. From \eqref{ineg1} we achieve
$$
\ln \E(\exp(tS(f)))\le Kn\overline\sigma^2_{k^\ast}(f)t^2,\mbox{ for }0\le t\le {k^\ast}^{-1},
$$
with $K=4(3e-4)$. Follow the Chernoff's device, i.e. using $\ln \P(S(f)\ge x)\le  \ln \E(\exp(tS(f)))-tx$ and optimizing in $0\le t\le {k^\ast}^{-1}$, we obtain
$$
 \P(S(f)\ge x)\le \exp\left(-\frac{x^2}{2Kn\overline\sigma^2_{k^\ast}(f)}\right)\1_{k^\ast x\le 2Kn\overline\sigma^2_{k^\ast}(f)} +\exp\left(\frac{Kn\overline\sigma^2_{k^\ast}(f)}{{k^\ast}^2}-\frac{x}{k^\ast}\right)\1_{k^\ast x> 2Kn\overline\sigma^2_{k^\ast}(f)}.
$$
Easy calculation yields  for all $x\ge0$
$$
\P(S(f)\ge\sqrt{2Kn\overline\sigma^2_{k^\ast}(f)}\1_{{k^\ast}^2 x\le 2Kn\overline\sigma^2_{k^\ast}(f)}+(k^\ast t+{k^\ast}^{-1}Kn\sigma^2_{k^\ast}(f))\1_{{k^\ast}^2 x> 2Kn\overline\sigma^2_{k^\ast}(f)}\le e{-x}.
$$
A rough bound $k^\ast t+{k^\ast}^{-1}Kn\sigma^2_{k^\ast}(f)\le 3k^\ast x/2$ for ${k^\ast}^2 x> 2Kn\overline\sigma^2_{k^\ast}(f)$ leads to the result of the Theorem.

\subsection{Proof of Proposition \ref{varest}}\label{pfvarest}
We have the classical decomposition
$$
 \vv\left(\sum_{i=1}^kf(X_i)\right)= k \vv(f(X_1))+2\sum_{r=1}^{k-1}(k-r)\cov(f(X_1),f(X_{r+1})).
$$
Now let us consider the coupling scheme $f(X_{r+1})^\ast$ distributed as $f(X_{r+1})$ but independent of $\mathcal M_1$. Then from Holder inequality it holds 
$$
\cov(f(X_1),f(X_{r+1}))=\E(\E(f(X_{r+1})-f(X_{r+1})^\ast\,|\,\mathcal M_1)f(X_1)).
$$
But as $f(X_{r+1})-f(X_{r+1})^\ast\le\delta_r$ conditionally to $\mathcal M_0$ we get the desired result.

\subsection{Proof of Proposition \ref{cb}}\label{pfcb}
We adapt the proof of \cite{Rio1996}. We are interested in estimated the coefficients $\varphi(\mathcal M_j,(X_{r+j},\ldots,X_{2r-1+j}))$ for any $(j,r)$ satisfying $1\le j\le j+r\le 2r-1+j\le n$. Let us fix $(j,r)$ and denote $(\xi^k_t)$ a sequence such that $\xi^k_{t}=\xi_{t}$ for all $t\ge r+j-k>j$ and $\xi_t^k=\xi_t'$ otherwise. Denote $U_t^k=F(\xi_{t-j}^k;j\in\N)$ and $X_t^k=H(U_t^k)$. For any $f\in\mathcal F$, we have
\begin{equation}\label{in0}
f(X_{r+j},\ldots,X_{2r-1+j})-f(X_{r+j}^k,\ldots,X_{2r-1+j}^k)\le \left(\frac{1}{r}\sum_{i=r+j}^{2r-1+j}d(X_i,X_i^k)\right)\wedge 1.
\end{equation}
By definition of the modulus of continuity and as $d(U_i^k,U_i)\le v_k$ for any $r+j\le i\le 2r-1+j$, we have 
$$
d(X_i,X_i^k)=d(H(U_i),H(U_i^k))\le w_H(U_i^k,v_k).
$$ 
Remarking that $\Big(r^{-1}\sum_{i=r+j}^{2r-1+j}w_H(U_i^k,v_k)\Big)\wedge 1$ is a measurable function of $((\xi'_t)_{t<r+j-k}, (\xi_t)_{t\ge r+j-k})$ bounded by $1$, it holds from the definition of the $\phi$-mixing coefficients:
$$
\E\left(\Big(r^{-1}\sum_{i=r+j}^{2r-1+j}w_H(U_i^k,v_k)\Big)\wedge 1~/~\mathcal M_j\right)\le \phi_{r-k}+\E\left(\Big(r^{-1}\sum_{i=r+j}^{2r-1+j}w_H(U_i^k,v_k)\Big)\wedge 1\right).
$$
Using again that $d(U_i^k,U_i)\le v_k$, then $w_H(U_i^k,v_k)\le 2w_H(U_i,2v_k)$. By stationarity of $(U_t)$, we obtain 
$$
\E\left(\Big(r^{-1}\sum_{i=r+j}^{2r-1+j}w_H(U_i^k,v_k)\Big)\wedge 1\right)\le \E(2w_H(U_0,2v_k))\wedge 1.
$$
So combining these inequalities we obtain for all $1\le k\le r-1$:
\begin{equation}\label{in1}
\left\|\E\left(f(X_{r+j},\ldots,X_{2r-1+j})-f(X_{r+j}^k,\ldots,X_{2r-1+j}^k)~\Big/~\mathcal M_j\right)\right\|_\infty\le\phi_{r-k}+ \E(2w_H(U_0,2v_k))\wedge 1.
\end{equation}
Using again the definition of the $\phi$-mixing coefficients, as $f$ is bounded by $1$ it holds
\begin{equation}\label{in2}
\left\|\E\left(f(X_{r+j}^k,\ldots,X_{2r-1+j}^k)~\Big/~\mathcal M_j\right)-\E\left(f(X_{r+j}^k,\ldots,X_{2r-1+j}^k)\right)\right\|_\infty\le\phi_{r-k}.
\end{equation}
Finally, using again \eqref{in0} and that $d(X_i,X_i^k)\le w_H(U_i,v_k)$, by stationarity of $(U_t)$ we obtain
\begin{equation}\label{in3}
\E f(X_{r+j},\ldots,X_{2r-1+j})- \E f(X_{r+j}^k,\ldots,X_{2r-1+j}^k)\le \E(w_H(U_0,v_k))\wedge 1.
\end{equation}
The result of the Proposition \ref{cb} follow from the definition of the $\varphi$-coefficients, the inequalities \eqref{in1}, \eqref{in2} and \eqref{in3}.

\subsection{Proof of Theorem \ref{thcl}}\label{pfthcl}
Let us keep the same notation than in the proof of Theorem \ref{mainsum}. The Benett's type deviation inequality follows classically from the Chernoff device applies with the estimate:
\begin{equation}\label{lapin}
\ln(\E(\exp(tS(f)))\le \frac{2n\sigma_k^2(f)}{k^2}(\exp(kt)-kt-1)+n\delta_k't\quad\mbox{for all }t\ge 0.
\end{equation}
To prove \eqref{lapin}, let us use the $L^\infty$-coupling scheme and \eqref{cl} to derive for all $1\le m\le p$:
$$ 
\left\|\sum_{i\in I_{2m}}f(X_i)-\sum_{i\in I_{2m}}f(X_i^\ast)\right\|_\infty\le \sum_{i\in I_{2(m+1)}}\left\|d(X_i,X_i^\ast)\right\|_\infty\le k\delta_k',
$$
where, as in Subsection \ref{pfmainsum}, $| I_j|=k$ for all $1\le j\le 2p$ with $2p-1\le nk^{-1}\le 2p$. Then, for all $t\ge 0$ we have:
$$
\exp\left(2t\sum_{i\in I_{2m}}f(X_i)\right)\le e^{2tk\delta_k'} \exp\left(2t\sum_{i\in I_{2m}}f(X_i^\ast)\right)\quad\mbox{a.s.}
$$
for all $1\le m\le p$. In particular, by independence of $(X_i^\ast)_{i\in I_{2m}}$ with $\mathcal M_{2{i-1}}$ and by stationary we deduce that
$$
\E\left(\exp\left(2t\sum_{i\in I_{2m}}f(X_i)\right)~|~\mathcal M_{2(m-1)}\right)\le e^{2tk\delta_k'} \E\left(\exp\left(2t\sum_{i\in I_{1}}f(X_i^\ast)\right)\right)
$$
for all $1\le m\le p$. Applying this inequality for $m=p$ we have
\begin{eqnarray*}
\E\exp(2tS_1)&=&\E\left(\exp\left(2t\sum_{1\le m\le p-1}\sum_{i\in I_{2m}}f(X_i)\right)\E\left(\exp\left(\sum_{i\in I_{2p}}f(X_i)\right)~|~\mathcal M_{2(p-1)}\right)\right)\\
&\le& e^{2tk\delta_k'}\E\left(\exp\left(2t\sum_{i\in I_{1}}f(X_i^\ast)\right)\right)\E\left(\exp\left(2t\sum_{1\le m\le p-1}\sum_{i\in I_{2m}}f(X_i)\right)\right).
\end{eqnarray*}
Let us do the same reasoning recursively on $m=p-1,\ldots,2$ to obtain finally
$$
\ln \E\exp(2tS_1)\le 2(p-1)k\delta_k't+p\ln\E\left(\exp\left(2t\sum_{i\in I_{1}}f(X_i^\ast)\right)\right).
$$
The classical Bennett inequality gives
$$
\ln\E\left(\exp\left(2t\sum_{i\in I_{1}}f(X_i^\ast)\right)\right)\le \frac{4\sigma_k^2(f)}{k}(\exp(kt)-kt-1)
$$
and the inequality \eqref{lapin} follows remarking that $4pk^{-1}\le 2nk^{-2}$ and $2(p-1)k\le n$.\\

For the Bernstein's type inequality, we use \eqref{lapin}, the series expansion of the function $\exp(x)-x-1$ and that $k!\ge 23^{k-2}$ for $k\ge2$ to derive:
$$
\ln(\E(\exp(tS(f)))\le \frac{n\sigma_k^2(f)t^2}{1-(k/3)t}+n\delta_k't\quad\mbox{for all }t\ge 0.
$$
With the same notation than in \cite{Massart2006}, for $x\ge n\delta'_k$ the Chernoff device leads to:
$$
\P(S(f)\ge x)\le \exp\left(\frac{2n\sigma_k^2(f)}{(k/3)^2}h_1\left(\frac{(k/3)(x-n\delta_k')}{2n\sigma_k^2(f)}\right)\right),
$$
where $h_1(x)=1+x-\sqrt{1+2x}$ for all $x\ge 0$. Then for all $x\ge 0$ we have 
$$
\P(S(f)\ge x+n\delta_k')\le \exp\left(\frac{2n\sigma_k^2(f)}{(k/3)^2}h_1\left(\frac{(k/3)x}{2n\sigma_k^2(f)}\right)\right)
$$
and the desired result follows as $h_1^{-1}(x)=\sqrt{2x}+x$ for all $x\ge 0$.

\subsection*{Acknowledgments}
The author is grateful to J\'er\^ome Dedecker  for
his helpful comments.

\end{document}